\font\math=msbm10 scaled 1000

\newcommand{\cpd}{{\hbox{\math o}}}

\documentstyle[11pt]{article}
 at12pt

\newcommand{\qed}{\hbox{\rule[-2pt]{3pt}{6pt}}}

\newtheorem{Thm}{Theorem}[section]

\newtheorem{ex}{Example}[section]

\newtheorem{defi}{Definition}[section]
\newtheorem{lem}{Lemma}[section]

\newtheorem{cor}{Corollary}[section]

\newcommand{\qedh}{\hfill\qed}

\newcommand{\vv}{ \vspace{.3in}  }

\begin{document}

\newpage

\begin{center}
\vspace*{.9in}

{\Large {\rm 
STRUCTURE OF THE GROUP OF \\
AUTOMORPHISMS
  OF C$^{*}$-ALGEBRAS
}}
  \end{center}

\vspace{.5in}
\begin{center}
{\large K.KAWAMURA}
\footnote{ e-mail : kawamura@kurims.kyoto-u.ac.jp.}

\end{center}
\begin{center}
{\normalsize  Research Institute for Mathematical Sciences }

{\normalsize  Kyoto University, Kyoto 606-8502, Japan}
\end{center}

\vspace{.4in}

\begin{quotation}\noindent A{\scriptsize BSTRACT}:
 We obtain a kind of structure theorem
 for  the automorphism group ${\rm Aut}{\cal A}$ 
 of a unital C$^{*}$-algebra ${\cal A}$.
 According to it, ${\rm Aut}{\cal A}$
 can be regarded as a subgroup
 of the semi-direct product of direct product 
 group consisting of
 some family of projective unitary groups
 and  some permutation group on the spectrum 
 of ${\cal A}$. 


\end{quotation}

\vspace{.4in}

\section{Introduction}

A non-commutative generalization of the functional representation theorem
for commutative unital C$^\ast$-algebras was introduced in \cite{CMP94}. This
generalization
was established via a non-commutative Gelfand transform mapping an unital
C$^\ast$-algebra 
$\cal A$ to an algebra of functions (for some non-commutative product) on
the set of 
pure states of $\cal A$ viewed as a uniform K\"ahler bundle over the spectrum 
of $\cal A$ (See Sect. 3). The K\"ahler structure involved can be seen as a
geometrical
counterpart of Shultz' characterization \cite{Shultz} of the set of pure 
states of a unital C$^\ast$-algebra.

As a consequence, any statement about C$^\ast$-algebras can be translated
into an equivalent
statement in terms of uniform K\"ahler bundles. For example, the set of
$^\ast$-isomorphisms
between two C$^\ast$-algebras $\cal A$ and $\cal A'$ is in one-to-one
correspondence 
with the set of uniform K\"ahler isomorphisms between the uniform K\"ahler
bundles associated
with $\cal A$ and $\cal A'$ \cite{CMP94}.

We think that this correspondence between C$^\ast$-algebras and
K\"ahler geometry can be advantageously exploited to get new insights in
some problems occurring in C$^\ast$-algebras theory. Also, the non-commutative
structure on the space of functions on the set of pure states seems to be
related
with deformation quantization of Poisson manifolds and a better understanding
of this link might result in a fruitful interaction between these fields.

In this paper, by using non-commutative functional representation theorem, 
we study the structure of the group of automorphisms
 of C$^{*}$-algebras in terms  of geometry of uniform K\"{a}hler bundles.

The paper is organized as follows.
In section 2, we state  our  main theorem.
 In section 3, we review the theory
 of uniform K\"{a}hler bundle \cite{CMP94}.
 In section 4, we introduce the orbit spectrum
 of a C$^{*}$-algebra ${\cal A}$. It is the space
 of orbits in the spectrum of ${\cal A}$
 of  the group of automorphisms of ${\cal A}$.
 We decompose the uniform K\"{a}{hler bundle
 associated to  ${\cal A}$ by
 the orbit spectrum. In section 5,
 we prove the main theorem.

\section{Structure of the group of automorphisms}

We first state our  main theorem using the following notation: \\
Let ${\cal A}$
 be a unital C$^{*}$-algebra.
${\rm Aut}{\cal A}$ is 
 the group of $^*$-automorphisms of ${\cal A}$,
 $B$  is  the spectrum
 of ${\cal A}$
 defined as the set of all the equivalence classes
 of irreducible representations of ${\cal A}$, and 
 ${\cal P}$ is  the set of pure states 
 of ${\cal A}$.  
With respect to the weak$^{*}$ topology, 
 ${\cal P}$ is a uniform space 
\cite{CMP94}, \cite{Bourbaki}.
 Since ${\rm Aut}{\cal A}$
 acts on $B$ naturally, we define   
 the orbit space $\Lambda\equiv 
 B/${\rm Aut}${\cal A}$ 
 denoting the corresponding natural projection 
 by $p^{'}:B\to \Lambda$. 

\begin{Thm}\label{Thm:main}
There is an injective homomorphism 
\[ \pi:
{\rm Aut}{\cal A}\hookrightarrow PU({\cal P})\cpd_{\delta}S(B)^{\Lambda}\]
 where
\[ PU({\cal P})\equiv \prod_{b\in B}PU_{b},\]
\[ S(B)^{\Lambda}\equiv \{ \phi: B\to B  :  \phi \mbox{ 
is a bijection 
 such that  } p^{'}\circ \phi=p^{'}\}, \]
 $PU_{b}$ is the projective unitary
 group  on the representation space  of the  representative element  $b$
 in the spectrum  $B$ 
 and 
 $\delta$ is the right action of $S(B)^{\Lambda}$ on $PU({\cal P})$
 defined by
\[ \{u_{b}\}\delta_{\phi}\equiv \{u_{\phi^{-1}(b)}\}\]
 for $\{u_{b}\}\in PU({\cal P})$,
 $u_{b}\in PU_{b}$ and $\phi\in S(B)^{\Lambda}$.

The image of ${\rm Aut}{\cal A}$  under $\pi$ is
 given in terms of a faithful  action $\kappa$
 of 
$PU({\cal P})\cpd_{\delta}S(B)^{\Lambda}$
 on ${\cal P}$, 
 by
\[ \{ g\in PU({\cal P})\cpd_{\delta}S(B)^{\Lambda}
 : \kappa_{g} \mbox{ is acting as a uniform homeomorphism on } {\cal P} \}.\]
 \end{Thm}

By this theorem, we characterize 
 an element of image of ${\rm Aut}{\cal A}$
 under $\pi$ 
as an element of $PU({\cal P})\cpd_{\delta}S(B)^{\Lambda}$
 which is a uniform homeomorphism on ${\cal P}$.

For $\alpha\in {\rm Aut}{\cal A}$,
 define $\alpha[\pi]\equiv [\pi\circ \alpha^{-1}]$
 for $[\pi]\in B$ where $[\pi]$
 is an equivalence class of
 irreducible representations with the  representative element
 $\pi$.

\begin{cor}
In Theorem \ref{Thm:main}, the image of the  subgroup
 $\{\alpha\in {\rm Aut}{\cal A}: 
 \alpha b=b \mbox{ for any } b\in B\}$ by $\pi$ is
 \[PU_{u}({\cal P})\equiv \{v\in PU({\cal P}):
 \kappa_{v} \mbox{ is a uniform homeomorphism on }
 {\cal P}\}\]
 where 
$PU_{u}({\cal P})$ is identified with  $PU_{u}({\cal P})\times\{ 1\}
\subset PU({\cal P})\cpd_{\delta}S(B)^{\Lambda}$.
\end{cor}

\begin{ex}{\rm
Let $X$ be a compact Hausdorff space. For the commutative C$^{*}$-algebra
 ${\cal A}\equiv C(X)$, ${\rm Aut}{\cal A}$
 is isomorphic to the group ${\rm Homeo}X$ of 
 homeomorphisms on $X$. 
So, 
 $\Lambda\equiv X/{\rm Homeo}X$
 depends on the topological structure of $X$. 
 Since 
 ${\cal P}\cong X\cong B$, 
$PU({\cal P}$
 is trivial. 
 So, $PU({\cal P})\cpd S(B)^{\Lambda}\cong S(B)^{\Lambda}$.

 Any compact Hausdorff space has uniformity and
 ${\rm Homeo}X$ is equal to the set of 
 uniform homeomorphisms on $X$ \cite{Bourbaki}.
 The image of the injection of ${\rm Aut}{\cal A}$
 into  $S(B_{{\cal A}})^{\Lambda}$
  is then equal to ${\rm Homeo}X$.

}
\end{ex}

By above argument, an element of 
$S(B_{{\cal A}})^{\Lambda}$ is considered as a "
 {\it topological} " symmetry of a 
 general noncommutative C$^{*}$-algebra ${\cal A}$.


%
%
\begin{ex}{\rm
For a C$^{*}$-algebra ${\cal A}$, ${\cal I}$
 is a {\it primitive } ideal of ${\cal A}$
 if there is an irreducible representation
 $\pi$
 of ${\cal A}$
 such that ${\rm ker}\pi={\cal I}$.
 The {\it primitive spectrum}  of ${\cal A}$
 is the set of all primitive spectrums of ${\cal A}$.
Assume ${\cal A}$ is a simple C$^{*}$-algebra.
 Then the primitive spectrum of ${\cal A}$
 consists 
of only one point. 
By definition,  the Jacobson
 topology of the spectrum $B$ is the trivial topology, 
 that is,  the 
 open sets  of $B$ are the empty set and  $B$ itself \cite{Ped}.
 So ${\rm Homeo}B=S(B)\equiv\{ \mbox{ permutation of $B$}\}$.
 Furthermore, if $\Lambda$ is 1-point ( we call ${\cal A}$  automorphic),
 then $S(B)^{\Lambda}=S(B)={\rm Homeo}B$.
 Thus, 
$PU({\cal P})\cpd S(B)^{\Lambda}
=PU({\cal P})\cpd {\rm Homeo}B$.
 Therefore ${\rm Aut}{\cal A}$ is a subgroup of 
$PU({\cal P})\cpd {\rm Homeo}B$
 if ${\cal A}$ is simple and automorphic. 
 
}
\end{ex}


%
%
\begin{ex}{\rm
For a C$^{*}$-algebra ${\cal A}$, ${\cal I}$
 is a {\it primitive } ideal of ${\cal A}$
 if there is an irreducible representation
 $\pi$
 of ${\cal A}$
 such that ${\rm ker}\pi={\cal I}$.
 The {\it primitive spectrum}  of ${\cal A}$
 is the set of all primitive spectrums of ${\cal A}$.
Assume ${\cal A}$ is a simple C$^{*}$-algebra.
 Then the primitive spectrum of ${\cal A}$
 consists 
of only one point. 
By definition,  the Jacobson
 topology of the spectrum $B$ is the trivial topology, 
 that is,  the 
 open sets  of $B$ are the empty set and  $B$ itself \cite{Ped}.
 So ${\rm Homeo}B=S(B)\equiv\{ \mbox{ permutation of $B$}\}$.
 Furthermore, if $\Lambda$ is 1-point ( we call ${\cal A}$  automorphic),
 then $S(B)^{\Lambda}=S(B)={\rm Homeo}B$.
 Thus, 
$PU({\cal P})\cpd S(B)^{\Lambda}
=PU({\cal P})\cpd {\rm Homeo}B$.
 Therefore ${\rm Aut}{\cal A}$ is a subgroup of 
$PU({\cal P})\cpd {\rm Homeo}B$
 if ${\cal A}$ is simple and automorphic. 
 
}
\end{ex}

\section{C$^{*}$-geometry}

 In this section, we review the characterization of 
 the set of pure states and the spectrum of a
 C$^{*}$-algebra following \cite{CMP94}.

Let $(f, E,M)$ be a surjective map $f:E\to M$ between
 two sets $E,M$.
\begin{defi}$(f, E,M)$
 is a formal K\"{a}hler bundle
 if 
there is a family $\{E_{m}\}_{m\in M}$ of K\"{a}hler manifolds
 indexed by $M$ and $E=\cup_{m\in M}E_{m}$ and $f(x)=m$ if $x\in E_{m}$.
\end{defi}

We simply denote $(f, E,M)$ by $E$.

Assume now that $E$ and $M$ are topological spaces.
\begin{defi}$(f, E,M)$
 is called a uniform K\"{a}hler bundle
 if $(f, E,M)$ is a formal K\"{a}hler bundle,
 $f$ is open,continuous, the topology of $E$
 is a uniform topology and the 
relative topology of each fiber is equivalent to 
 the K\"{a}hler topology of its fiber.
\end{defi}

For a uniform topology, see \cite{Bourbaki}. 
The weak$^{*}$-topology
 on the set of pure states of C$^{*}$-algebra is 
   a uniform topology.

\begin{defi}
Two formal K\"{a}hler bundle $(f, E,M)$, $(f^{'}, E^{'},M^{'})$ are
 isomorphic
 if there is  a pair $(\beta,\phi)$
 of bijections $\beta :E\to E^{'}$ and $\phi:M\to M^{'}$,
 such that 
 $f^{'}\circ\beta=\phi\circ f$
\[ 
\begin{array}{ccccc}
&& \beta &&\\
&E&\cong & E^{'}&\\
f & \downarrow   & 
& \downarrow & f^{'}\\
&M & \cong & M^{'}&  \\
&& \phi &&
\end{array}
\]
and any 
restriction $\beta |_{f^{-1}(m)} :f^{-1}(m)\to (f^{'})^{-1}(\phi(m))$
 is a holomorphic K\"{a}hler isometry for any $m\in M$.
We call $(\beta,\phi)$ a
 formal K\"{a}hler isomorphism between  $(f, E,M)$ and $(f^{'}, E^{'},M^{'})$.
\end{defi}

By definition of a formal K\"{a}hler bundle
 isomorphism $(\beta,\phi)$
 between $(f, E,M)$ and $(f^{'}, E^{'},M^{'})$, $\phi$
 is uniquely determined by  $\beta$:
For $m\in M$, the value $\phi(m)$
 is given by $\phi(m)=f^{'}(\beta(e))$
 with arbitrary $e\in f^{-1}(\{m\})$.

\begin{defi}
Two uniform K\"{a}hler bundles $(f, E,M)$, $(f^{'}, E^{'},M^{'})$ are
 isomorphic
 if there is  a formal K\"{a}hler isomorphism
 $(\beta,\phi)$ between  $(f, E,M)$ and $(f^{'}, E^{'},M^{'})$
 such that $\phi$ is a homeomorphisms,  
 and $\beta$ is a uniform homeomorphism.
We call $(\beta,\phi)$ a
 uniform K\"{a}hler isomorphism between  
$(f, E,M)$ and $(f^{'}, E^{'},M^{'})$. 
\end{defi}

By definition, any uniform K\"{a}hler bundle is a 
 formal K\"{a}hler bundle. For 
 a uniform K\"{a}hler bundle $E$, we define:
\begin{defi}\label{defi:iso1}
\[\widetilde{{\rm Iso}(E)}\equiv \mbox{
 the group of formal K\"{a}hler bundle isomorphisms of $E$}, \]
\[ {\rm Iso}(E) \equiv\mbox{ 
  the group of uniform  K\"{a}hler bundle isomorphisms of $E$} .\]
\end{defi}

 By the GNS representation, there is a natural
 projection $p:{\cal P}\to B$
 from the set ${\cal P}$ of pure states onto
 the spectrum $B$. We consider $(p,{\cal P},B)$
  as  a map of topological spaces where
 ${\cal P}$ is endowed  with weak$^{*}$ topology and $B$ is endowed with 
 the Jacobson topology.

In Ref.\cite{CMP94}, the following results are proved.
\begin{Thm}( Reduced atomic realization )
 For any unital C$^{*}$-algebra ${\cal A}$,
$(p,{\cal P},B)$ is  a
 uniform K\"{a}hler bundle.
\end{Thm}
 
For a fiber ${\cal P}_{b}\equiv p^{-1}(b)$, let $(\pi_{b},{\cal H}_{b})$
 be some irreducible representation
belonging to $b\in B$. To $\rho\in {\cal P}_{b}$, correspond 
 $[x_{\rho}]\in {\cal P}({\cal H}_{b})\equiv 
({\cal H}_{b}\setminus\{0\})/{\bf C}^{\times}$
 where $\rho=\omega_{x_{\rho}}\circ \pi_{b}$
 with $\omega_{x_{\rho}}$ denoting a vector state
$\omega_{x_{\rho}}=<x_{\rho}|(\cdot)x_{\rho}>$.
 Then ${\cal P}_{b}$ has a K\"{a}hler manifold structure 
 induced by this correspondence from projective Hilbert space
 ${\cal P}({\cal H}_{b})$.

\begin{Thm}
Let ${\cal A}_{i}$ be C$^{*}$-algebras
 with associated uniform K\"{a}hler bundles 
$(p_{i}, {\cal P}_{i}, B_{i})$, $i=1,2$.
 Then 
${\cal A}_{1}$ and ${\cal A}_{2}$ are $^*$-isomorphic 
 if and only if  $(p_{1}, {\cal P}_{1}, B_{1})$
 and  $(p_{2}, {\cal P}_{2}, B_{2})$
 are isomorphic as uniform K\"{a}hler bundles.
\end{Thm}

\begin{cor}\label{cor:iso1}
Let ${\rm Aut}{\cal A}$ be the group of $^*$-automorphisms
 of a C$^{*}$-algebra ${\cal A}$ with an associated
 uniform K\"{a}hler bundle
 ${\cal P}=(p, {\cal P},B)$,
 and ${\rm Iso}{\cal P}$ be the group of uniform K\"{a}hler bundle 
 automorphisms on ${\cal P}$.
 Then
 there is a group isomorphism
\[ {\rm Aut}{\cal A}\cong {\rm Iso}{\cal P}.\]

\end{cor}

For $\alpha\in {\rm Aut}{\cal A}$, let 
 $\beta_{\alpha}\equiv \alpha^{*}|_{{\cal P}} :{\cal P}\to {\cal P}$,
 $\alpha^{*}(\rho)\equiv \rho\circ\alpha^{-1}$
 and induced bijection $\phi_{\alpha}:B\to B$
 defined  by
 $\phi_{\alpha}([\pi])\equiv [\pi\circ \alpha^{-1}]$.
 Then $(\beta_{\alpha},\phi_{\alpha})$ 
 becomes a uniform K\"{a}hler bundle automorphism of ${\cal P}$.

We call these objects {\it $C^{*}$-geometry}
 since any C$^{*}$-algebra can be 
 reconstructed from the associated uniform K\"{a}hler bundle
 \cite{CMP94} and, therefore,
 any C$^{*}$-algebra is determined
 by such a geometry.

By the above result, we can consider the structure
 of ${\rm Aut}{\cal A}$ in the language  of ${\rm Iso}{\cal P}$.

\section{Orbit decomposition of a K\"{a}hler bundle}

We decompose the set of pure states and the spectrum of a C$^{*}$-algebra
 ${\cal A}$ as a uniform K\"{a}hler bundle. 
By using
 this decomposition, we describe automorphisms of ${\cal A}$
 in each decomposed component in the next section.

In \cite{K-R}(II, p 906),  two pure states $\rho$ and $\rho^{'}$ of 
 a C$^{*}$-algebra ${\cal A}$ are  called {\it automorphic}
 if  there is an automorphism $\alpha$ of ${\cal A}$ 
 such that $\rho^{'}=\rho\circ\alpha$. For example,
 any two pure states of a uniform matricial algebra
 are automorphic ( \cite{K-R} II, Theorem 12.3.4 ).
 In general, the set of pure states
 of a C$^{*}$-algebra
 is divided into a disjoint union of automorphic component.
 Therefore, each automorphism induces transformations 
 on  each automorphic components.
 The idea on which  this section 
 is based comes from this point of view.

 Let $G\equiv{\rm Aut}{\cal A}$.
 $G$ is naturally acting on the spectrum $B$
 by $g[\pi]\equiv[\pi\circ g^{-1}]$
 for $g\in G$ and $[\pi]\in B$. 
So,   
 we define the space $\Lambda$ of orbits of $G$ in $B$,
\[ \Lambda\equiv B/G,\]
 and call it the {\it orbit spectrum}.
$G$ acts naturally also  on ${\cal P}$ by $g\rho\equiv \rho\circ g^{-1}$.
%
%
\begin{lem}\label{lem:onetoone}
 The orbit of $G$ in ${\cal P}$ and $B$ are 
 in one-to-one correspondence.
\end{lem}
{\it Proof.}
  Let $G\rho$ be an orbit of $G$ through $\rho\in {\cal P}$.
 We define 
$\Psi(G\rho)\equiv G[\pi_{\rho}]\in \Lambda$ where
 $[\pi_{\rho}]$ is the unitary equivalence class
 of irreducible representations of ${\cal A}$
 with a representative element $\pi_{\rho}$, 
 given by the GNS representation of $\rho$.
For $g\rho\in G\rho$,
 $\pi_{g\rho}$ is  unitarily equivalent to 
 $g\pi_{\rho}\equiv\pi_{\rho}\circ g^{-1}\in G[\pi_{\rho}]$
 by uniqueness of the GNS representation.
 Then the map $\Psi: {\cal P}/G\to \Lambda$
 is well defined. By definition, $\Psi$
 maps orbits in ${\cal P}$ to orbits in $B$. 
And 
$\Psi(G\rho)=Gp(\rho)$. 
Hence, $\Psi$ is onto.

 If $\Psi(G\rho)=\Psi(G\rho^{'})$ and, 
 $(\pi_{\rho},{\cal H}_{\rho},x_{\rho})$ and 
$(\pi_{\rho^{'}},{\cal H}_{\rho^{'}},x_{\rho^{'}})$
 are GNS representations of $\rho$, $\rho^{'}\in {\cal P}$
 respectively,
  then, $G[\pi_{\rho}]=G[\pi_{\rho^{'}}]$.
Since two automorphic pure states have
  GNS representation spaces with the same dimension, 
 there are $g\in G$, a  representative element
 $\pi^{'}\in [\pi_{\rho^{'}}]$
 which acts on ${\cal H}_{\rho}$,
 $\rho^{'}=\omega_{x_{\rho}}\circ\pi^{'}$ and 
 a unitary operator $U$ on  ${\cal H}_{\rho}$
 such that $\pi^{'}={\rm Ad}U\circ g\pi_{\rho}$.
By irreducibility of $\pi_{\rho}$,
 we can choose 
 a unitary element $V$ in ${\cal A}$ such that 
\[\rho^{'}=(g\rho)\circ {\rm Ad}V=({\rm Ad}V^{*}\circ g)\rho\in G\rho\]
 ( see \cite{K-R}, II, 10.2.6.).
 Therefore $G\rho=G\rho^{'}$. $\Psi$ is an injection. 
\[
\begin{array}{ccccc}
{\cal P} & \stackrel{p}{\to} & B & \stackrel{p^{'}}{\to} & B/G=\Lambda \\
 &&&& \\
\downarrow && & \nearrow _{\Psi}&\\
&&&& \\
{\cal P}/G &&&&
\end{array}.
\]
\qedh

We decompose ${\cal P}$
 by $\Lambda$ into the family of  uniform K\"{a}hler bundles.

Let $p^{'}:B\to\Lambda=B/{\rm Aut}{\cal A}$ be the natural projection
 with   fibers 
 given by $B^{\lambda}\equiv(p^{'})^{-1}(\lambda)$, 
 $\lambda\in \Lambda$. 
 Let ${\cal P}^{\lambda}\equiv \cup_{b\in B^{\lambda}}{\cal P}_{b}$.
 By  lemma \ref{lem:onetoone}, $B^{\lambda}$ and ${\cal P}^{\lambda}$
 are orbits of $G$ in  $B$ and ${\cal P}$,  respectively.
Let $p^{\lambda}\equiv p|_{{\cal P}^{\lambda}}$ for $\lambda\in\Lambda$.
 Then $(p^{\lambda}, {\cal P}^{\lambda}, B^{\lambda})$ 
  for each $\lambda\in \Lambda$
 becomes a uniform K\"{a}hler bundle with the relative topology
such that its total space ${\cal P}^{\lambda}$ 
is automorphic, that is,
 any
 two elements of ${\cal P}^{\lambda}$ are
 transformed by some automorphisms  of ${\cal A}$ each other.
We obtain a decomposition
\[(p,{\cal P},B)=\bigcup_{\lambda\in \Lambda}(p^{\lambda},{\cal P}^{\lambda},B^{\lambda})\]
 of a uniform K\"{a}hler bundle.

 Any two elements $b,b^{'}$ in the same 
 orbit $B^{\lambda}$ have representative
 representation spaces with the same dimensions.
 For an orbit $\lambda\in \Lambda$,
  let ${\cal H}_{\lambda}$ be a Hilbert space
  corresponding to a representative element
 of some point in an orbit $B^{\lambda}$.
We can choose  a representative element  belonging to $B^{\lambda}$
 which acts on the same Hilbert space ${\cal H}_{\lambda}$. 

Let ${\cal P}^{o}\equiv \cup_{\lambda\in \Lambda}
{\cal P}({\cal H}_{\lambda})\times B^{\lambda}$ and 
$p^{o}:{\cal P}^{o}\to B$ defined
 by $p^{o}(\xi,b)=b$ for $\xi\in {\cal P}({\cal H}_{\lambda})$
 and $b\in B^{\lambda}$. $(p^{o},{\cal P}^{o}, B)$
 becomes a formal K\"{a}hler bundle
 with fiber ${\cal P}({\cal H}_{\lambda})\times\{b\}$
 for $b\in B^{\lambda}\subset B$.
\begin{Thm}
$(p,{\cal P},B)$ and  $(p^{o},{\cal P}^{o}, B)$
 are isomorphic as formal K\"{a}hler bundles.
\end{Thm}
{\it Proof.}
Fix a family of representative elements $\{\pi_{b}\}_{b\in B}$ of $B$
 such that $\pi_{b}$ acts on ${\cal H}_{\lambda}$
 if $b\in B^{\lambda}$.
Define  $\beta:{\cal P}\to {\cal P}^{0}$
 by $\beta(\rho)\equiv([x_{\rho}],b)
\in {\cal P}({\cal H}_{\lambda})\times B^{\lambda}$
 if $[\pi_{\rho}]=b\in B^{\lambda}$, 
where $\rho=\omega_{x_{\rho}}\circ \pi_{b}$
 and $x_{\rho}\in {\cal H}_{\lambda}$.
 Let $\phi:B\to B$ be the identity map on $B$.
 Then $\beta({\cal P}_{b})={\cal P}({\cal H}_{\lambda})\times\{b\}$
 for $b\in B^{\lambda}$. And $\beta$ is fiber-wise holomorphic isometry.
Then, $(\beta,\phi)$ becomes a formal K\"{a}hler isomorphism
 between $(p,{\cal P},B)$ and  $(p^{o},{\cal P}^{o}, B)$. \qedh
 
\begin{cor}\label{cor:obdeco}(Orbit decomposition)
Let ${\cal P}=(p,{\cal P},B)$ be
 a uniform K\"{a}hler bundle associated with a C$^{*}$-algebra.
 Then there is a uniform K\"{a}hler bundle ${\cal P}^{o}=
(p^{o},{\cal P}^{o}, B)$
 with ${\cal P}^{o}=\cup_{\lambda\in \Lambda}
{\cal P}({\cal H}_{\lambda})\times B^{\lambda}$ such that
 ${\cal P}\cong {\cal P}^{o}$.
\end{cor}
{\it Proof.}
In the previous theorem, let the topology of  ${\cal P}^{o}$
 be the induced topology from ${\cal P}$. Then $(p^{o},{\cal P}^{o}, B)$
 becomes a uniform K\"{a}hler bundle which
 is isomorphic to $(p,{\cal P},B)$. \qedh

\section{Proof of the main theorem}

By Corollary \ref{cor:obdeco}, we identify the uniform K\"{a}hler bundle
 ${\cal P}=(p,{\cal P},B)$ associated with a C$^{*}$-algebra
 ${\cal A}$ and its orbit decomposition
 ${\cal P}^{o}$ corresponding to the orbit
 spectrum $\Lambda$.
Let ${\rm Aut}{\cal A}$ be the group
 of $^*$-automorphisms of a C$^{*}$-algebra ${\cal A}$
 with associated uniform K\"{a}hler bundle ${\cal P}=(p,{\cal P},B)$
 and the orbit spectrum $\Lambda$. 
 By corollary \ref{cor:obdeco}, we identify ${\cal P}$
 with  its orbit decomposition ${\cal P}^{o}$.

 Recall that $PU({\cal P})$, $S(B)^{\Lambda}$
 and  $\widetilde{{\rm Iso}{\cal P}}$ 
are defined in 
Theorem \ref{Thm:main}
 and Definition \ref{defi:iso1}.

We define actions  $t,s$ of $PU({\cal P})$, 
$S(B)^{\Lambda}$ respectively on the set ${\cal P}$
 of pure states  by
\[
\begin{array}{ccccccc}
t_{u}(\xi,b)    &    \equiv & ( & u_{b}\xi &,   &b & ), \\
s_{\phi}(\xi,b) &     \equiv & ( & \xi&,    &\phi(b) & ), \\
\end{array}
\]
 for $u=\{u_{b}\}_{b\in B}\in PU({\cal P})$,
 $\phi\in S(B)^{\Lambda}$ 
and $(\xi,b)\in {\cal P}^{o}$.
With these actions, we define injective homomorphisms
 $\tau$ and $\sigma$ of $PU({\cal P})$ and of $S(B)^{\Lambda}$
 into $\widetilde{{\rm Iso}{\cal P}}$,
 by
\[ 
\begin{array}{ccccccc}
\tau_{u} & \equiv & ( & t_{u} &,& id_{B} & ),\\
\sigma_{\phi} & \equiv & ( & s_{\phi} & , & \phi & ), 
\end{array}
\]
 respectively, for $u\in PU({\cal P})$,
 $\phi\in S(B)^{\Lambda}$.
\begin{lem}\label{lem:iso33}
 \[\tau(PU({\cal P}))=\{ (\beta,\phi)\in \widetilde{{\rm Iso}{\cal P}}:
\phi=id_{B}
\}\equiv G_{3}.\]
\end{lem}
{\it Proof.}
By definition
 of $\tau$,  $\tau(PU({\cal P}))$ is contained in
 $G_{3}$. On the other hand,  for any $g=(\beta, id_{B})\in G_{3}$, 
$\beta$ 
 becomes a holomorphic  K\"{a}hler isometry
 on each fiber
 by definition of $G_{3}$. So it becomes a projective unitary
 on each fiber
 by Wigner's theorem. 
Thus, there is a family of
 projective unitaries  corresponding  to $g$.  \qedh


For $u\in PU({\cal P})$
 and $\phi\in S(B)^{\Lambda}$,
 we obtain $s_{\phi}t_{u}s_{\phi^{-1}}=t_{u\delta_{\phi}}$,
 where $\delta$ is the right action of $S(B)^{\Lambda}$
 on $PU({\cal P})$
 defined in Theorem \ref{Thm:main}.
 From this follows the relation
\begin{equation}\label{equation:cova}
\sigma_{\phi}\tau_{u}\sigma_{\phi^{-1}}=\tau_{u\delta_{\phi}}.
\end{equation}

Consider the  action $\tilde{\delta}_{a}={\rm Ad}_{a}$ of
 $a\in \sigma (S(B)^{\Lambda})$ on $\tau( PU({\cal P}))$.
 Let $G_{2}$ be the  group generated by 
$\tau( PU({\cal P}))$ and $\sigma(S(B)^{\Lambda})$.
 We obtain the following isomorphism between semi-direct products
 of groups.
\begin{lem}\label{lem:embed1}
\[   PU({\cal P})\cpd_{\delta} S(B)^{\Lambda}
\cong \tau( PU({\cal P}))\cpd_{\tilde{\delta}}\sigma(S(B)^{\Lambda})
=G_{2}.\]
\end{lem}
{\it Proof.}
By \ref{equation:cova} and the 
 definition of semi-direct product, 
 the lemma follows.
\qedh

{\it Proof of Theorem \ref{Thm:main}( main Theorem ).}
By lemma \ref{lem:embed1}, $PU({\cal P})\cpd_{\delta} S(B)^{\Lambda}$
 is embedded into $\widetilde{{\rm Iso}{\cal P}}$
 as a subgroup.

Let 
\begin{equation}\label{eqn:kappa}
\kappa:PU({\cal P})\cpd_{\delta} S(B)^{\Lambda}
 \hookrightarrow \widetilde{{\rm Iso}{\cal P}}
\end{equation}
 be defined by $\kappa(u,\phi)\equiv\tau_{u}\sigma_{\phi}
=(t_{u}s_{\phi},\phi)$ for
 $(u,\phi)\in PU({\cal P})\cpd_{\delta} S(B)^{\Lambda}$.

On the other hand, by Corollary \ref{cor:iso1}, there is
 an  isomorphism
\begin{equation}\label{eqn:pi1}
 \pi_{1}:{\rm Aut}{\cal A} \cong{\rm Iso}{\cal P}\subset  
\widetilde{{\rm Iso}{\cal P}}.
\end{equation}
 We denote
$\pi_{1}(\alpha)\equiv (\beta_{\alpha},\phi_{\alpha})$
 for $\alpha\in {\rm Aut}{\cal A}$. Then,
 $(\beta_{\alpha}\circ s_{\phi_{\alpha}^{-1}},id_{B})\in G_{3}$.
 By lemma \ref{lem:iso33}, there is $u^{\alpha}\in PU({\cal P})$
 such that
$\beta_{\alpha}\circ s_{\phi_{\alpha}^{-1}}=t_{u^{\alpha}}$.
So, 
$\beta_{\alpha}=t_{u^{\alpha}}s_{\phi_{\alpha}}$.
Therefore, 
\[\pi_{1}(\alpha)= (\beta_{\alpha},\phi_{\alpha})
=(t_{u^{\alpha}}s_{\phi_{\alpha}}, \phi_{\alpha})
=\tau_{u^{\alpha}}\sigma_{\phi_{\alpha}}.\]
By this calculation
 and \ref{eqn:kappa}, $\pi_{1}({\rm Aut}{\cal A})\subset G_{2}
=\kappa(PU({\cal P})\cpd_{\delta} S(B)^{\Lambda})$.

 Let $\pi\equiv \kappa^{-1}|_{G_{2}}\circ \pi_{1}$.
 We denote an element of 
$PU({\cal P})\cpd_{\delta} S(B)^{\Lambda}$
 by $u\cdot\delta_{\phi}$
 which satisfies the product law
\[ (u\cdot\delta_{\phi})(u^{'}\cdot\delta_{\phi^{'}})
=\{u(u^{'}\delta_{\phi})\}\cdot\delta_{\phi\phi^{'}}\]
 for $u\cdot\delta_{\phi}$,   $u^{'}\cdot\delta_{\phi^{'}}
\in PU({\cal P})\cpd_{\delta} S(B)^{\Lambda})$,
 $u,u^{'}\in PU({\cal P})$
 and $\phi,\phi^{'}\in S(B)^{\Lambda}$.
 Then 
\[\pi(\alpha)=u^{\alpha}\cdot \delta_{\phi^{\alpha}}\in
 PU({\cal P})\cpd_{\delta} S(B)^{\Lambda} \]
 for $\alpha\in {\rm Aut}{\cal A}$.
 Then we obtain the  injective homomorphism 
 for which we have been looking 
 \[
\begin{array}{ccccc}
& &\pi && \\
{\rm Aut}{\cal A}&&\hookrightarrow && PU({\cal P})\cpd_{\delta} S(B)^{\Lambda}. \\
&&&&\\
\pi_{1}&\searrow & & \nearrow& \kappa^{-1}|_{G_{2}}\\
&&&&\\
&&G_{2}&&
\end{array}
\]
By Corollary \ref{cor:iso1}, we obtain
 \[\pi({\rm Aut}{\cal A})= (\kappa^{-1}|_{G_{2}}\circ \pi_{1})({\rm Aut}{\cal A})
=\kappa^{-1}|_{G_{2}}({\rm Iso}{\cal P})\mbox{ (by equation \ref{eqn:pi1}) } \]
\[=\kappa^{-1}|_{G_{2}}(\{ (\beta,\phi)\in \widetilde{{\rm Iso}{\cal P}}:
 \beta \mbox{ is a uniform  homeomorphism on } {\cal P} \} )\]
\[=\{ g\in  PU({\cal P})\cpd_{\delta} S(B)^{\Lambda}:
 \kappa_{g}\mbox{ acts on ${\cal P}$
 as a uniform  homeomorphism}\},\]
from which the statement
 of Theorem \ref{Thm:main}
 immediately follows.
\qedh

\section{Conclusion}

In this paper, we have obtained  the orbit decomposition
 of the uniform K\"{a}hler bundles
 and the group of automorphisms.

The next step would be  to consider
  the orbit decomposition of the algebra itself.
 In this context, the meaning of decomposition 
 has to be cleared.
 It might  be a decomposition like by crossed products,
 free products of C$^{*}$-algebras.

We are studying geometrical objects corresponding to 
 modules, crossed product, subalgebra, $^*$-homomorphism
 and etc  are currently under study for non-commutative
 C$^{*}$-algebras. 
They are realization of non-commutative geometry
 by `` {\it real }" geometry
 defined as the
 set of points and its function space
\cite{CMP94}. So, they must be direct generalization
 of the geometry of commutative case
 leading to  new geometrical structures 
for which we would like to give a better understanding.

{\bf Acknoledgement}
 I thank I.Ojima and G.Dito for a critical reading of the paper.

\vv

\newpage

\end{document}